\documentclass{article}
\usepackage[a4paper, total={7in, 10in}]{geometry}
\usepackage[T1]{fontenc}
\usepackage{graphicx}
\usepackage{mathtools}
\usepackage{amsmath,amsthm,amssymb,mathrsfs,amstext, titlesec,enumitem, comment, graphicx, color, xcolor, stmaryrd,mathabx}
\usepackage{thmtools}
\usepackage{xpatch}
\usepackage{tikz-cd}
\usepackage{nameref}
\usepackage{hyperref}
\makeatletter
\def\Ddots{\mathinner{\mkern1mu\raise\p@
\vbox{\kern7\p@\hbox{.}}\mkern2mu
\raise4\p@\hbox{.}\mkern2mu\raise7\p@\hbox{.}\mkern1mu}}
\makeatother
\newtheorem{theorem}{Theorem}[section]
\newtheorem{corollary}[theorem]{Corollary}
\newtheorem{lemma}[theorem]{Lemma}
\newtheorem{question}[theorem]{Question}
\theoremstyle{definition}
\newtheorem{definition}[theorem]{Definition}

\begin{document}
	
	\title{ A combinatorial approach to exponential patterns in multiplicative $IP^{\star}$ sets in $\mathbb{N}$}

	\date{}
	\author{Pintu Debnath
		\footnote{Department of Mathematics,
			Basirhat College,
			Basirhat-743412, North 24th parganas, West Bengal, India.\hfill\break
			{\tt pintumath1989@gmail.com}}
	}
	\maketitle

\begin{abstract}

In [On $IP^{\star}$sets and central sets, Combinatorica, 14 (1994)  269-277], N. Hindman and V.Bergelson proved additive $IP^{\star}$-sets contain finite sums and finite products of a single sequence. An analogous study was made by A. Sisto in [Exponential triples, Electronics Journal of Combinatorics, 18 (2011), no. 147], where he proved that multiplicative $IP^{\star}$-sets contain exponential $IP$ of type $I$ and finite sums of a single sequence as well as exponential $IP$ of type $II$ and finite products of another single sequence,
 using the algebra in the Stone-\v{C}ech Compactification of discrete semigroups. In this article, we will provide a combinatorial proof of the result of A. Sisto.

\end{abstract}

 

\maketitle

\section{Introduction}
 At the beginning of this section, we start with some definitions to better explain what we are going to do in this article. 
\begin{definition} Let $A$ be a subset of $\mathbb{N}$.
\begin{itemize}
   \item[(a)] (\textbf{Additive $IP$ set}): The set $A$ is called additive $IP$ set if for a sequence $\langle x_{n}\rangle_{n=1}^{\infty}$ in $\mathbb{N}$ such thst $FS\left(\langle x_{n}\rangle_{n=1}^{\infty}\right)\subseteq A$, where $FS\left(\langle x_{n}\rangle_{n=1}^{\infty}\right)=\left\{ \sum_{t\in H}x_{t}:H\in\mathcal{P}_{f}\left(\mathbb{N}\right)\right\}$.
   
   \item[(b)]  (\textbf{Multiplicative $IP$ set}): The set $A$ is called additive $IP$ set if for a sequence $\langle x_{n}\rangle_{n=1}^{\infty}$ in $\mathbb{N}$ such that  $FP\left(\langle x_{n}\rangle_{n=1}^{\infty}\right)\subseteq A$, where $FP\left(\langle x_{n}\rangle_{n=1}^{\infty}\right)=\left\{ \prod_{t\in H}x_{t}:H\in\mathcal{P}_{f}\left(\mathbb{N}\right)\right\}$.
   
   \item[(c)] (\textbf{Additive $IP^{\star}$ set}): The set $A$ is called additive $IP^{\star}$ set if it intersects with all additive $IP$ set.
   
   \item[(d)] (\textbf{Multiplicative $IP^{\star}$ set}): The set $A$ is called multiplicative $IP^{\star}$ set if it intersects with all multiplicative $IP$ set.

\end{itemize}
\end{definition}

$IP^{\star}$ sets play a crucial and fundamental role in arithmetic Ramsey theory and number theory. Some recent results  \cite{G1},\cite{GHW}, on prime number theory  are directly connected with additive  $IP^{\star}$ set.

\begin{theorem}{(\textbf{Finite Sums Theorem})}
 Let $r\in \mathbb{N}$ and let $\mathbb{N}=C_{1}\cup C_{2}\cup\ldots\cup C_{r}$.  There exist $i\in \left\{1,2,\ldots,r\right\}$ and a sequence  $\langle x_{n}\rangle_{n=1}^{\infty}$ in $\mathbb{N}$ such that $FS\left(\langle x_{n}\rangle_{n=1}^{\infty}\right)\subseteq C_{i}$.
    
\end{theorem}

The Finite Sums Theorem was a significantly strong conjecture of Graham and Rothschild and   proved by N.Hindman \cite{H} using   combinatorial method in 1974. In the same year J.E.Baumgarthner, in \cite{B} provided a short  combinatorial proof of the Hindman Finite Sums Theorem. Using algebra of $\beta\mathbb{N}$, an algebraic proof of the  Hindman Finite Sums Theorem is followed from \cite[Corollary 5.10 Page-112]{HS}, but it was originaly due to Glazer and Galvin.

The Hindman Finite Sum Theorem simply state that any finite coloring of the set of positive integers $\mathbb{N}$, there is a monocromatic copy of additive $IP$ set. Passing to the map $n\rightarrow{2^{n}}$ for each $n\in\mathbb{N}$, we immediately have a monocromatic copy of multiplicative $IP$ set. In \cite[Theorem 2.4]{BH93}, N.Hindman and V.Bergelson combinatorially proved the following and algebraic proof is given in \cite[Corollary 5.22 Page-116]{HS}:

\begin{theorem}\cite[Theorem 2.4]{BH93}
 Let $r\in \mathbb{N}$ and let $\mathbb{N}=C_{1}\cup C_{2}\cup\ldots\cup C_{r}$.  There exist $i\in \left\{1,2,\ldots,r\right\}$ and two sequences  $\langle x_{n}\rangle_{n=1}^{\infty}$ and $\langle y_{n}\rangle_{n=1}^{\infty}$ in $\mathbb{N}$ such that $FS\left(\langle x_{n}\rangle_{n=1}^{\infty}\right)\cup FP\left(\langle y_{n}\rangle_{n=1}^{\infty}\right)\subseteq C_{i}$.
    
\end{theorem}

 Natural question arises in our mind that whether we get a single sequence in the conclusion of the above theorem. From \cite[Theorem 17.16]{HS}, we get a negative answer.

 Following definition from \cite[Definition 5.13 Page-112]{HS}.
 
\begin{definition} 
Let $\langle x_{n}\rangle_{n=1}^{\infty}$ be a sequence in $\mathbb{N}$.
\begin{itemize}
    \item[(a)] (\textbf{Sum subsystem}): The sequence $\langle y_{n}\rangle_{n=1}^{\infty}$ is called sum
subsystem of $\langle x_{n}\rangle_{n=1}^{\infty}$ if and only if
there is a sequence $\langle H_{n}\rangle_{n=1}^{\infty}$ in $\mathcal{P}_{f}\left(\mathbb{N}\right)$
such that for every $n\in\mathbb{N}$, $\max H_{n}<\min H_{n+1}$  and $y_{n}=\sum_{t\in H_{n}}x_{t}$.
    \item[(b)] (\textbf{Product subsystem}): The sequence $\langle y_{n}\rangle_{n=1}^{\infty}$ is called product
subsystem of $\langle x_{n}\rangle_{n=1}^{\infty}$ if and only if
there is a sequence $\langle H_{n}\rangle_{n=1}^{\infty}$ in $\mathcal{P}_{f}\left(\mathbb{N}\right)$
such that for every $n\in\mathbb{N}$, $\max H_{n}<\min H_{n+1}$  and $y_{n}=\prod_{t\in H_{n}}x_{t}$.
\end{itemize}

\end{definition}
 The following theorem due to Bergelson and Hindman \cite[Theorem 2.6]{BH94} by combinatorial method and algebraic proof follows from \cite[Corollary 16.21 Page-413]{HS}.
\begin{theorem}\label{AIP* finite sum product}\cite[theorem 2.6]{BH94}
Let $\langle x_{n}\rangle_{n=1}^{\infty}$ be a sequence in $\mathbb{N}$
and $A$ be an additive $IP^{*}$ set in $\left(\mathbb{N},+\right)$. Then there
exists a sum subsystem $\langle y_{n}\rangle_{n=1}^{\infty}$ of $\langle x_{n}\rangle_{n=1}^{\infty}$
such that $FS\left(\langle y_{n}\rangle_{n=1}^{\infty}\right)\cup FP\left(\langle y_{n}\rangle_{n=1}^{\infty}\right)\subseteq A.$
\end{theorem}

As we mentioned in the abstract that A. Sisto did an analogous study for multiplicative $IP^{\star}$ sets in $\mathbb{N}$, we need the following definition to state the work of A. Sisto.

\begin{definition}\label{exponential definition}\cite[Definition 2]{S} Let $\langle x_{n}\rangle_{n=1}^{\infty}$ be sequence in $\mathbb{N}$. Define inductively

\begin{itemize}
    \item[(a)]  $FE^{I}\left(\langle x_{n}\rangle_{n=1}^{k+1}\right)=\left\{y^{x_{k+1}}|y\in FE^{I}\left(\langle x_{n}\rangle_{n=1}^{k}\right) \right\}\cup FE^{I}\left(\langle x_{n}\rangle_{n=1}^{k}\right)\cup \left\{x_{k+1}\right\} $,
    \item[(b)]  $FE^{II}\left(\langle x_{n}\rangle_{n=1}^{k+1}\right)=\left\{\left(x_{k+1}\right)^{y}|y\in FE^{II}\left(\langle x_{n}\rangle_{n=1}^{k+1}\right) \right\}\cup FE^{II}\left(\langle x_{n}\rangle_{n=1}^{k}\right)\cup \left\{x_{k+1}\right\} $,
\end{itemize}
with $FE^{I}\left(\langle x_{n}\rangle_{n=1}^{1}\right)=FE^{I}\left(\langle x_{n}\rangle_{n=1}^{1}\right)=\left\{x_{1}\right\}$. Set
\begin{itemize}
    \item[(c)] $FE^{I}\left(\langle x_{n}\rangle_{n=1}^{\infty}\right)=\bigcup_{k\in \mathbb{N}}FE^{I}\left(\langle x_{n}\rangle_{n=1}^{k}\right)$ and
    \item[(d)] $FE^{II}\left(\langle x_{n}\rangle_{n=1}^{\infty}\right)=\bigcup_{k\in \mathbb{N}}FE^{II}\left(\langle x_{n}\rangle_{n=1}^{k}\right).$
\end{itemize}
We will say that $C\subseteq\mathbb{N}$ is an exponential $IP$-set of type $I$(resp. $II$) if it contains a set $FE^{I}\left(\langle x_{n}\rangle_{n=1}^{\infty}\right)$(resp. $FE^{II}\left(\langle x_{n}\rangle_{n=1}^{\infty}\right)$) for some sequence $\langle x_{n}\rangle_{n=1}^{\infty}$.
\end{definition}
 In \cite[Theorem 3]{S},  A. Sisto proved the following theorem algebraically.
 \begin{theorem}\label{sisto's result}
     Given any multiplicative $IP^{\star}$ set $A$,  there exist some sequences $\langle x_{n}\rangle_{n=1}^{\infty}$ and  $\langle y_{n}\rangle_{n=1}^{\infty}$ in  $\mathbb{N}$ such that
     \begin{itemize}
         \item[(a)] $FS\left(\langle x_{n}\rangle_{n=1}^{\infty}\right)\cup FE^{I}\left(\langle x_{n}\rangle_{n=1}^{\infty}\right)\subseteq A$ and
         \item[(b)]$ FP\left(\langle y_{n}\rangle_{n=1}^{\infty}\right)\cup FE^{II}\left(\langle y_{n}\rangle_{n=1}^{\infty}\right)\subseteq A$.
     \end{itemize}
 \end{theorem}
In \textbf{section 2}, we will provide a simple combinatorial proof of Theorem \ref{AIP* finite sum product} and by the same technique we will give a combinatorial proof of Theorem\ref{sisto's result} in \textbf{section 3} .

\section{Finite sums and finite products in additive $IP^{\star}$ sets}

In \cite[Theorem 2]{B}, James E. Baumgartaner gave a short proof of the Hindman Finite Sums Theorem combinatorially. Actually, James E.Baumgartaner established the following theorem, from which we get the Hindman Finite Sums Theorem. 
\begin{theorem}\cite[Theorem 2]{B}
Let $r\in\mathbb{N}$ and $\mathcal{P}_{f}\left(\mathbb{N}\right)=C_{1}\cup C_{2}\cup\ldots\cup C_{r}$
be a finite partition of $\mathcal{P}_{f}\left(\mathbb{N}\right)$.
Then there exist $i\in\left\{ 1,2,\ldots,r\right\} $ and a sequence
  $\langle K_{n}\rangle_{n=1}^{\infty}$  in $\mathcal{P}_{f}\left(\mathbb{N}\right)$ such that $FU\left(\langle K_{i}\rangle_{n=1}^{\infty}\right)\subseteq C_{i}$
with $K_{i}\cap K_{j}=\emptyset$ if $i\neq j$.
\end{theorem}

Using the above theorem, we get the following trivial corollary:
\begin{corollary}\label{partition of set}
Let $r\in\mathbb{N}$ and $\mathcal{P}_{f}\left(\mathbb{N}\right)=C_{1}\cup C_{2}\cup\ldots\cup C_{r}$
be a finite partition of $\mathcal{P}_{f}\left(\mathbb{N}\right)$.
Then there exist $i\in\left\{ 1,2,\ldots,r\right\} $ and a sequence
$\langle H_{n}\rangle_{n=1}^{\infty}$ in $\mathcal{P}_{f}\left(\mathbb{N}\right)$ such that $FU\left(\langle H_{i}\rangle_{n=1}^{\infty}\right)\subseteq C_{i}$
with $\max H_{i}<\min H_{i+1}$.
\end{corollary}

The following is a very well known result. We present the proof as we have no specific reference to the elementary proof of the following:

\begin{corollary}\label{partition of IP set}
Let $r\in\mathbb{N}$ and $FS\left(\langle x_{n}\rangle_{n=1}^{\infty}\right)=C_{1}\cup C_{2}\cup\ldots\cup C_{r}$
be a finite partition of $\mathbb{N}$. Then there exist $i\in\left\{ 1,2,\ldots,r\right\} $
and a sequence $\langle y_{n}\rangle_{n=1}^{\infty}$ of sum subsystem
of $\langle x_{n}\rangle_{n=1}^{\infty}$ such that $FS\left(\langle y_{n}\rangle_{n=1}^{\infty}\right)\subseteq C_{i}$.
\end{corollary}

\begin{proof}
We may consider that $\langle x_{n}\rangle_{n=1}^{\infty}$ is an injective sequence. Let  $\langle z_{n}\rangle_{n=1}^{\infty}$ be a subsequence of $\langle x_{n}\rangle_{n=1}^{\infty}$
with the property $z_{1}+z_{2}+\cdots+z_{n}<z_{n+1}$,  for all $n\in\mathbb{N}$.
 Then we get,
 \begin{itemize}
     \item[(1)]  $\sum_{n\in A}z_{n}=\sum_{n\in B}z_{n}\iff A=B$ for all $n\in\mathbb{N}$ and
     \item[(2)]  $FS\left(\langle z_{n}\rangle_{n=1}^{\infty}\right)=C_{1}^{\prime}\cup C_{2}^{\prime}\cup\ldots\cup C_{r}^{\prime}$, where $C_{i}^{\prime}=C_{i}\cap FS\left(\langle z_{n}\rangle_{n=1}^{\infty}\right)$ for all $i\left\{\in\{1,2,\ldots,r\right\}$.
 \end{itemize}

Construct a partition $\mathcal{P}_{f}\left(\mathbb{N}\right)=C_{1}^{\prime\prime}\cup C_{2}^{\prime\prime}\cup\ldots\cup C_{r}^{\prime\prime}$
by $A\in C_{i}^{\prime\prime}\iff\sum_{t\in A}z_{t}\in C_{i}^{\prime}$.
 By the  Corollary \ref{partition of set}, we get $i\in\left\{ 1,2,\ldots,r\right\} $
and a sequence $\langle H_{n}\rangle_{n=1}^{\infty}$  in $\mathcal{P}_{f}\left(\mathbb{N}\right)$ such that $FU\langle H_{n}\rangle_{n=1}^{\infty}\subseteq C_{i}^{\prime\prime}$
with $\max H_{n}<\min H_{n+1}$ for all $n\in \mathbb{N}$. From the construction of partition
of $\mathcal{P}_{f}\left(\mathbb{N}\right)$, we get a sum subsystem $\langle z_{n}^{\prime}\rangle_{n=1}^{\infty}$
of $\langle z_{n}\rangle_{n=1}^{\infty}$ such that
$FS\langle z_{n}^{\prime}\rangle_{n=1}^{\infty}\subseteq C_{i}^{\prime}\subseteq C_{i}$  for some $i\in\left\{1,2,\ldots,r\right\}$ . As $\langle z_{n}\rangle_{n=1}^{\infty}$ is a subsequence of $\langle x_{n}\rangle_{n=1}^{\infty}$, then $\langle z_{n}^{\prime}\rangle_{n=1}^{\infty}$ is a sum subsystem of $\langle x_{n}\rangle_{n=1}^{\infty}$. To get the result take $z^{\prime}_{n}=y_{n}$ for all $n\in \mathbb{N}$.
 
\end{proof}

An additive  $IP^{\star}$ set meets with every additive $IP$ set and a much stronger statement is the following:
\begin{corollary}\label{intere section AIP* and AIP}
Let $\langle x_{n}\rangle_{n=1}^{\infty}$ be a sequence in $\mathbb{N}$
and $A$ be an additive $IP^{*}$ set in $\left(\mathbb{N},+\right)$. Then there exists a sum subsystem $\langle y_{n}\rangle_{n=1}^{\infty}$ of $\langle x_{n}\rangle_{n=1}^{\infty}$ such that $FS\left(\langle y_{n}\rangle_{n=1}^{\infty}\right)\subseteq A$.
\end{corollary}

\begin{proof}
Let $FS\langle x_{n}\rangle_{n=1}^{\infty}=C_{1}\cup C_{2}$, where
$C_{1}=A\cap FS\left(\langle x_{n}\rangle_{n=1}^{\infty}\right)$
and $C_{2}=FS\left(\langle x_{n}\rangle_{n=1}^{\infty}\right)\setminus A$.
Use the corollary \ref{partition of IP set} and the fact that $A$ is an additive  $IP^{\star}$ set.
\end{proof}
\begin{corollary}\label{intersection of IP* set}
Let $A,B\subseteq\mathbb{N}$ are two additive $IP^{*}$ sets in $\mathbb{N}$ then $A\cap B$ is an  additive $IP^{*}$ set in $\mathbb{N}$.
\end{corollary}

\begin{proof}
For any sequence $\langle x_{n}\rangle_{n=1}^{\infty}$, by Corollary \ref{intere section AIP* and AIP},  we get that there exists a sum subsystem $\langle y_{n}\rangle_{n=1}^{\infty}$
of $\langle x_{n}\rangle_{n=1}^{\infty}$ such that $FS\left(\langle y_{n}\rangle_{n=1}^{\infty}\right)\subseteq A$.
In the same argument, we get a sum subsystem $\langle z_{n}\rangle_{n=1}^{\infty}$
of $\langle y_{n}\rangle_{n=1}^{\infty}$ such that $FS\left(\langle z_{n}\rangle_{n=1}^{\infty}\right)\subseteq B$
and $FS\left(\langle z_{n}\rangle_{n=1}^{\infty}\right)\subseteq A\cap B$.
\end{proof}
Now, we are in a position to prove the main result of this section.

\begin{definition}
    Let $n\in \mathbb{N}$. If $A\subseteq\mathbb{N}$ denote
    \begin{itemize}
        \item[(a)] $-n+A=\left\{m\in\mathbb{N}:m+n\in A\right\}$,
        \item[(b)] $n^{-1}A=\left\{m\in\mathbb{N}:mn\in A\right\}$,  
    \end{itemize}
\end{definition}

Let $A$ be an additive $IP^{\star}$ set and $y\in\mathbb{N}$. Then by \cite[Lemma 16.19 Page-411]{HS},  $y^{-1}A$ is an additive $IP^{\star}$ set in $\mathbb{N}$.

\begin{proof}[\textbf{Proof of Theorem \ref{AIP* finite sum product}}.]
 By Corollary \ref{intere section AIP* and AIP},  for any sequence $\langle x_{n}\rangle_{n=1}^{\infty}$ we get a sum subsystem
$\langle z_{n}\rangle_{n=1}^{\infty}$ of $\langle x_{n}\rangle_{n=1}^{\infty}$
such that $FS\left(\langle z_{n}\rangle_{n=1}^{\infty}\right)\subseteq A$.
Here $z_{n}=\sum_{t\in H_{n}}x_{t}$ with $\max H_{n}<\min H_{n+1}$ for all $n\in \mathbb{N}$.
Let $B=FS\left(\langle z_{n}\rangle_{n=1}^{\infty}\right)$ and $y_{1}=z_{1}$. Then
$FS\left(\langle z_{i}\rangle_{i=2}^{\infty}\right)\subseteq B\cap\left(-y_{1}+B\right)$. Pick $y_{2}\in FS\left(\langle z_{i}\rangle_{i=2}^{\infty}\right)\cap y_{1}^{-1}A$
with $y_{2}=\sum_{t\in K_{2}}z_{t}$ and $K_{2}\in\mathcal{P}_{f}\left(\mathbb{N}\right)$. Hence $\left\{y_{1},y_{2},y_{1}+y_{2},y_{1}y_{2}\right\}\subset A$. Inductively, let $m\in\mathbb{N} $ and assume that we have chosen $\langle y_{n}\rangle_{n=1}^{m}$ and $\langle K_{n}\rangle_{n=1}^{m}$ such that:

\begin{itemize}
    \item[(1)] each $y_{n}=\sum_{t\in H_{n}}z_{t}$,
    \item[(2)] if $n<m$, then $\max K_{n}<\min K_{n+1}$, and
    \item[(3)] $FS\left(\langle y_{n}\rangle_{n=1}^{m}\right)\cup FE\left(\langle y_{n}\rangle_{n=1}^{m}\right)\subset A$.
 \end{itemize}

Let $C=FS\left(\langle y_{n}\rangle_{n=1}^{m}\right)$ and $D= FP\left(\langle y_{n}\rangle_{n=1}^{m}\right)$. Then
\begin{itemize}
    \item[(4)] $FS\left(\langle z_{i}\rangle_{i=\max K_{m}+1}^{\infty}\right)\subseteq B\cap\bigcap_{y\in C}\left(-y+B\right)$ and pick
    \item[(5)] $y_{m+1}\in FS\left(\langle z_{i}\rangle_{i=\max K_{m}+1}^{\infty}\right)\cap\bigcap_{y\in D}y^{-1}A.$

\end{itemize}

Then $y_{m+1}=\sum_{t\in K_{m+1}}z_{t}$ with $\max K_{m} < \min K_{m+1}$ and $$\left\{ y_{m+1} \right\}\cup \left\{y+ y_{m+1}:y\in FS\left(\langle y_{n}\rangle_{n=1}^{m}\right)  \right\}\cup\left\{y y_{m+1}:y\in FP\left(\langle y_{n}\rangle_{n=1}^{m}\right) \right\}\subset A.$$

  By the method of  induction, we get a sequence $\langle y_{n}\rangle_{n=1}^{\infty}$ in $\mathbb{N}$ such that  $FS\left(\langle y_{n}\rangle_{n=1}^{\infty}\right)\cup FP\left(\langle y_{n}\rangle_{n=1}^{\infty}\right)\subseteq A$, 
 where $y_{n}=\sum_{t\in K_{n}}z_{t}=\sum_{t\in K_{n}}\sum_{s\in H_{t}}x_{s}=\sum_{t\in G_{n}}x_{t}$ and 
 $G_{n}=\cup_{s\in K_{n}}H_{s}$ with $\max G_{n}<\min G_{n+1}$ for all $n\in \mathbb{N}$.
\end{proof}

\section{Combined exponential patterns in  multiplicative $IP^{\star}$ sets}
 In the previous section, we have discussed on combined additive and multiplicative patterns in additive $IP^{\star}$ set. In this section we discussed on Sisto's work regarding multiplicative $IP^{\star}$ set  in \cite{S}.
Now, start with the following definition:
\begin{definition}
    If $A\subseteq\mathbb{N}$ denote
    \begin{itemize}
        \item[(a)] if $n\geq 2$, $\log_{n}[A]=\left\{m\in\mathbb{N}:n^{m}\in A\right\}$,
        \item[(b)] if $n\geq 1$, $A^{1/n}=\left\{m\in\mathbb{N}:n^{m}\in A\right\}$,
        
    \end{itemize}
\end{definition}
\begin{lemma}\cite[Lemma 13]{S}
    Let $A$ be a multiplicative $IP^{\star}$-set and let $n\in\mathbb{N}$.
    \begin{itemize}
        \item[(a)] If $n\geq 2$ then $\log_{n}[A]$ is an additive $IP^{\star}$-set.
        \item[(b)] If $n\geq 1$ then $A^{1/n}$ is a multiplicative $IP^{\star}$-set.
    \end{itemize}
\end{lemma}
The following simple result is essential to proving the  main results in this section:
\begin{lemma}\label{MIP* implies AIP}
    Let $A$ be a  multiplicative $IP^{\star}$ set. Then $A$ is an additive $IP$ set.
\end{lemma}
\begin{proof}
    We first prove that $A$ is multiplicative syndetic set. If possible let Let $A$ is not multiplicative syndetic. Pick $x_{1}  \in  \mathbb{N}\setminus A$ and $x_{2} \in  \mathbb{N}\setminus A\cup x_{1}^{-1}A$, then $\left\{x_{1},x_{2},x_{1}x_{2}\right\}\subset A$. By induction, we can find out a sequence $\langle x_{n}\rangle _{n=1}^{\infty}$ in $\mathbb{N}$, such that
$FP\left(\langle x_{n}\rangle_{n=1}^{\infty}\right)\cap A=\emptyset$, where
 $x_{n}  \in  \mathbb{N}\setminus A\cup\bigcup_{x\in FP\left(\langle x_{k}\rangle_{k=1}^{n-1}\right)}x^{-1}A $
for all $n\in\mathbb{N}$. Which contradicts the fact that $A$ is $MIP^{\star}$-set.

As $A$ is multiplicative syndetic, there exists a finite sequence $\langle y_{n}\rangle_{n=1}^{r}$ such that $\mathbb{N}=y_{1}^{-1}A\cap y_{2}^{-1}A\cap\ldots\cap y_{r}^{-1}A $. Then by Hindman Finite Sums Theorem, there exist a sequence $\langle z_{n}\rangle_{n=1}^{\infty}$ and  $k\in\left\{1,2,\ldots,r\right\}$, such that $FS\left(\langle z_{n}\rangle_{n=1}^{\infty}\right)\subseteq y_{k}^{-1}A$, which implies $FS\left(\langle y_{k}z_{n}\rangle_{n=1}^{\infty}\right)\subseteq A$
\end{proof}

\begin{theorem}
     Given any multiplicative $IP^{\star}$ set $A$ there exists some infinite sequence $\langle y_{n}\rangle_{n=1}^{\infty}$  in  $\mathbb{N}$ such that $FS\left(\langle y_{n}\rangle_{n=1}^{\infty}\right)\cup FE^{I}\left(\langle y_{n}\rangle_{n=1}^{\infty}\right)\subseteq A$.
 \end{theorem}
 
 \begin{proof}
     As per Lemma\ref{MIP* implies AIP}, there exists a sequence $\langle x_{n}\rangle_{n=1}^{\infty}$ in $\mathbb{N}$ such that $FS\left(\langle x_{n}\rangle_{n=1}^{\infty}\right)\subseteq A$. Let $B=FS\left(\langle x_{n}\rangle_{n=1}^{\infty}\right)$ and $y_{1}=x_{1}$. Then
$FS\left(\langle x_{i}\rangle_{i=2}^{\infty}\right)\subseteq B\cap\left(-y_{1}+B\right)$. Pick $y_{2}\in FS\left(\langle x_{i}\rangle_{i=2}^{\infty}\right)\cap \log_{y_{1}}[A]$
with $y_{2}=\sum_{t\in K_{2}}x_{t}$ and $K_{2}\in\mathcal{P}_{f}\left(\mathbb{N}\right)$. Hence $\left\{y_{1},y_{2},y_{1}+y_{2},{y_{1}}^{y_{2}}\right\}\subset A$. Inductively, let $m\in\mathbb{N} $ and assume that we have chosen $\langle y_{n}\rangle_{n=1}^{m}$ and $\langle K_{n}\rangle_{n=1}^{m}$ such that:

\begin{itemize}
    \item[(1)] each $y_{n}=\sum_{t\in H_{n}}z_{t}$,
    \item[(2)] if $n<m$, then $\max K_{n}<\min K_{n+1}$, and
    \item[(3)] $FS\left(\langle y_{n}\rangle_{n=1}^{m}\right)\cup FE^{I}\left(\langle y_{n}\rangle_{n=1}^{m}\right)\subset A$.
 \end{itemize}

Let $C=FS\left(\langle y_{n}\rangle_{n=1}^{m}\right)$ and $D= FP\left(\langle y_{n}\rangle_{n=1}^{m}\right)$. Then
\begin{itemize}
    \item[(4)] $FS\left(\langle x_{i}\rangle_{i=\max K_{m}+1}^{\infty}\right)\subseteq B\cap\bigcap_{y\in C}\left(-y+B\right)$ and pick
    \item[(5)] $y_{m+1}\in FS\left(\langle z_{i}\rangle_{i=\max K_{m}+1}^{\infty}\right)\cap\bigcap_{y\in D}\log_{y}[A].$

\end{itemize}

Then $y_{m+1}=\sum_{t\in K_{m+1}}x_{t}$ with $\max K_{m} < \min K_{m+1}$ and $$\left\{ y_{m+1} \right\}\cup \left\{y+ y_{m+1}:y\in FS\left(\langle y_{n}\rangle_{n=1}^{m}\right)  \right\}\cup\left\{y^{ y_{m+1}}:y\in FE^{I}\left(\langle y_{n}\rangle_{n=1}^{m}\right) \right\}\subset A.$$

  By the method of  induction, we get a sequence $\langle y_{n}\rangle_{n=1}^{\infty}$ in $\mathbb{N}$ such that  $FS\left(\langle y_{n}\rangle_{n=1}^{\infty}\right)\cup FE^{I}\left(\langle y_{n}\rangle_{n=1}^{\infty}\right)\subseteq A$, 
 where $y_{n}=\sum_{x\in K_{n}}x_{t}$  with  $\max K_{n}<\min K_{n+1}$ for all $n\in \mathbb{N}$.
\end{proof}

In the Theorem \ref{AIP* finite sum product}, it is observed that the existing sequence $\langle y_{n}\rangle_{n=1}^{\infty}$ is a sum subsystem of the given sequence $\langle x_{n}\rangle_{n=1}^{\infty}$ but we are unable to find out this type of result for the above theorem. Although, the following theorem reflects the analogous result of the Theorem \ref{AIP* finite sum product}.

\begin{theorem}
    Let $\langle x_{n}\rangle_{n=1}^{\infty}$ be a seqience in $\mathbb{N}$. Given any multiplicative $IP^{\star}$ set $A$ there exists a product subsystem  $\langle y_{n}\rangle_{n=1}^{\infty}$  of $\langle x_{n}\rangle_{n=1}^{\infty}$  such that $FP\left(\langle y_{n}\rangle_{n=1}^{\infty}\right)\cup FE^{II}\left(\langle y_{n}\rangle_{n=1}^{\infty}\right)\subset A$.
 \end{theorem}
 
 \begin{proof}
     Let $\langle x_{n}\rangle_{n=1}^{\infty}$  be a sequence $\mathbb{N}$. There exists a product subsystem
$\langle z_{n}\rangle_{n=1}^{\infty}$ of $\langle x_{n}\rangle_{n=1}^{\infty}$
such that $FP\left(\langle z_{n}\rangle_{n=1}^{\infty}\right)\subseteq A$. By the same process in \ref{AIP* finite sum product}, replacing $-y+A$ by $y^{-1}A$  and replacing $y^{-1}A$ by $A^{1/y}$, we get the proof.
 \end{proof}
Combining the above two theorems, we get the Theorem \ref{sisto's result} and the  following question arises naturally:

 \begin{question}
     Given any multiplicative $IP^{\star}$ set $A$, does there exist a  sequences $\langle x_{n}\rangle_{n=1}^{\infty}$  in  $\mathbb{N}$ such that
     \begin{itemize}
         \item[(a)] $FS\left(\langle x_{n}\rangle_{n=1}^{\infty}\right)\cup FE^{I}\left(\langle x_{n}\rangle_{n=1}^{\infty}\right)\subseteq A$
         \item[(b)]$ FP\left(\langle x_{n}\rangle_{n=1}^{\infty}\right)\cup FE^{II}\left(\langle x_{n}\rangle_{n=1}^{\infty}\right)\subseteq A$?
     \end{itemize}
     
 \end{question}
 We do not know the answer of the above question. But in \cite[Theorem 1.4]{DG}, the author of this article and Goswami provided the following partial answer using algebra of $\beta\mathbb{N}$.
\begin{theorem}
Let $A$ be a multiplicative $IP^{\star}$ set. Then there exists a sequence $\langle x_{n}\rangle_{n=1}^{\infty}$
such that 
$$ FP\left(\langle x_{n}\rangle_{n=1}^{\infty}\right)\cup FE^{I}\left(\langle x_{n}\rangle_{n=1}^{\infty}\right)\cup FE^{II}\left(\langle x_{n}\rangle_{n=1}^{\infty}\right)\subseteq A.$$

\end{theorem}
As our aim in this article is to prove all results combinatorially. It is a natural question: can we prove the above theorem combinatorially? We will provide a partial answer and combinatorially prove the following:
\begin{theorem}\label{Partial answer}
Let $A$ be a multiplicative $IP^{\star}$ set. Then there exists a sequence   $\langle x_{n}\rangle_{n=1}^{\infty}$ in $\mathbb{N}$ 
such that $$FE^{I}\left(\langle x_{n}\rangle_{n=1}^{\infty}\right)\bigcup FE^{II}\left(\langle x_{n}\rangle_{n=1}^{\infty}\right)\subseteq A.$$
\end{theorem}

\begin{lemma}
    Let $A$ be an additive $IP^{\star}$ set and $B$ be an multiplicative $IP^{\star}$ set in $\mathbb{N}$. Then the intersection of $A$ and $B$ is nonempty.
\end{lemma}
\begin{proof}
    By  Lemma\ref{MIP* implies AIP}, $B$ is additive $IP$ set. Then the intersection of $A$ and $B$ is nonempty.
\end{proof}

\begin{proof}[\textbf{Proof of Theorem \ref{Partial answer}}.]
    Let $x_{1}\in A$. Then $\log_{x_{1}}[A]$ is an additive $IP^{\star}$-set and $A^{1/{x_{1}}}$ is a multiplicative $IP^{\star}$-set. As $A\cap A^{1/{x_{1}}} $ is multiplicative $IP^{\star}$-set, we have $A\bigcap A^{1/{x_{1}}}\bigcap \log_{x_{1}}[A]\neq\phi $. 
    Picking $x_{2}\in A\cap A^{1/{x_{1}}}\bigcap \log_{x_{1}}[A]\neq\phi$, we get $\left\{ x_{1},x_{2},x_{1}^{x_{2}},x_{2}^{x_{1}}\right\}\subset A$. 
    
    To prove by induction, let $m\in\mathbb{N}$ and finite sequence $\langle x_{n}\rangle_{n=1}^{m}$,  $FE^{I}\left(\langle x_{n}\rangle_{n=1}^{m}\right)\bigcup FE^{II}\left(\langle x_{n}\rangle_{n=1}^{m}\right)\subseteq A$. Then $A\cap\bigcap_{x\in  FE^{II}\left(\langle x_{n}\rangle_{n=1}^{m}\right)}A^{1/x}$ is multiplicative $IP^{\star}$-set and $\bigcap_{x\in  FE^{I}\left(\langle x_{n}\rangle_{n=1}^{m}\right)}\log_{x}[A]$ is additive $IP^{\star}$-set. So, there intersection is non empty, take 
\begin{itemize}
    \item $x_{m+1}\in A\cap\bigcap_{x\in  FE^{II}\left(\langle x_{n}\rangle_{n=1}^{m}\right)}A^{1/x}\bigcap_{x\in  FE^{I}\left(\langle x_{n}\rangle_{n=1}^{m}\right)}\log_{x}[A]$, which implies
    \item $\left\{x_{m+1}\right\}\cup \left\{x_{m+1}^{x}:x\in FE^{II}\left(\langle x_{n}\rangle_{n=1}^{m}\right) \right\}\cup \left\{x^{x_{m+1}}:x\in FE^{I}\left(\langle x_{n}\rangle_{n=1}^{m}\right)\right\}\subset A$ i.e., 
    \item $FE^{I}\left(\langle x_{n}\rangle_{n=1}^{m+1}\right)\bigcup FE^{II}\left(\langle x_{n}\rangle_{n=1}^{m+1}\right)\subseteq A$.
    
\end{itemize}
    
\end{proof}

We conclude this section with the following question:

\begin{question}
Let $A$ be a multiplicative $IP^{\star}$ set.  Does there exist a sequence $\langle x_{n}\rangle_{n=1}^{\infty}$
such that 
$$ FS\left(\langle x_{n}\rangle_{n=1}^{\infty}\right)\cup FE^{I}\left(\langle x_{n}\rangle_{n=1}^{\infty}\right)\cup FE^{II}\left(\langle x_{n}\rangle_{n=1}^{\infty}\right)\subseteq A ?$$
\end{question}

\end{document}